\documentclass[10pt,letterpaper,twocolumn]{CLS/autart}
\usepackage[latin9]{inputenc}
\usepackage{color}
\definecolor{note_fontcolor}{rgb}{0.800781, 0.800781, 0.800781}
\usepackage{float}
\usepackage{amsmath}
\usepackage{amssymb}
\usepackage{graphicx}

\makeatletter

\usepackage{alltt}
\usepackage{subfigure}
\usepackage{graphics}
\usepackage{epsfig}
\usepackage{psfrag}
\usepackage{makeidx}
\usepackage[british]{babel}

\usepackage[margin=2cm]{geometry}

\usepackage{amscd,amsfonts}
\usepackage{graphicx}
\usepackage{tikz}
\usetikzlibrary{decorations.pathreplacing,calc}
\usetikzlibrary{matrix,decorations.pathreplacing,calc}

\usepackage{color}
\newcommand{\fin}{\hfill$\Box$}

\newcommand{\n}{\normalfont}

\newtheorem{lemma}{Lemma}
\newtheorem{deff}{Definition}
\newtheorem{obs}{Remark}

\makeatother

\begin{document}
\vskip -5cm \begin{frontmatter}
\title{Moment matching based reduced closed-loop design to achieve  asymptotic performance\thanksref{footnoteinfo}}
\thanks[footnoteinfo]{The research  has received funding from: GAR2023 code 71 Research Grant funded and managed by the Patrimony Foundation (Fundatia "Patrimoniu") of the Romanian Academy, from the Recurrent Fund of Donors, Contract No.  260/28.11.2023;  National Program for 
Research of the National Association of Technical Universities - GNAC ARUT 2023 grant, contract no. 5/06.10.2023; UEFISCDI, Romania,  PN-III-P4-PCE-2021-0720, under project L2O-MOC, nr. 70/2022. This paper has not been presented at any IFAC meeting.}

%author[Kawano]{Yu Kawano}\ead{ykawano@hiroshima-u.ac.jp},
\author[Ionescu_UPB,Ionescu_ISMMA]{Tudor C. Ionescu}\ead{tudor.ionescu@acse.pub.ro}   %and 
%\author[Ionescu_UPB]{Cristian Oar\u a}\ead{cristian.oara@acse.pub.ro}

\address[Ionescu_UPB]{Department of Automatic Control and Systems Engineering, Politehnica Univ. of Bucharest, 060042 Bucharest, Romania.} 
\address[Ionescu_ISMMA]{"Gheorghe Mihoc-Caius Iacob" Institute of Mathematical Statistics and Applied Mathematics of the Romanian Academy, 050711 Bucharest, Romania.} 
%\address[Kawano]{Hiroshima University, Japan}

\begin{keyword}                            
Moment matching, controller, low-order,  closed-loop, reduced order, asymptotic tracking, performance. 
\end{keyword} 
                                                                       
\vspace{-0.5cm} 

\begin{abstract} 
In this paper, the moment matching techniques are adopted to obtain  reduced-order closed-loop systems with reduced-order  controllers that maintain the closed-loop stability and guarantee desired asymptotic performance, after revealing the relationship between the Internal Model Principle used in control design and the time-domain moment matching problem. As a result, the design of a low order controller can be done starting from considering the achieving of asymptotic performance as a moment matching problem, resulting in a reduced order closed-loop system.  
\end{abstract}
\end{frontmatter}

\section{Introduction}

\label{intro}

In the field of linear control design, a low-order controller is often
obtained by computing a high-order controller at first and then reducing
the order of the controller by using some techniques. The problem
of controller reduction has long been studied in the literature in
different model reduction frameworks, see e.g., \cite{sandberg-murray-OCAM2009}.
Most of the earlier results are based on stability preserving balanced
truncation, see e.g., \cite{sandberg-TAC2010}. Many algorithms result
in various controller representations, meeting different performance
requirements; see \cite{anderson-liu-TAC1987}. Variations of the
problem and the results have been introduced in, e.g., \cite{schelfourt-demoor-TAC2002},
where it has been shown that closed-loop balanced truncation with
preservation of the closed-loop properties is actually equivalent
to frequency weighted balanced truncation with certain selections
of the weighting systems. A different approach based on interpolation
approaches has been taken in, e.g., \cite{reis-styke-MCDMS2007,vandendorpe-vandooren-2009}.
A Krylov-based approach to controller reduction has been used in \cite{gugercin-antoulas-beattie-gildin-CDC2004},
where the reduced-order controller  is guaranteed to match the desired
closed-loop behavior, maintaining stability. More recently data-driven control strategies have been developed in the Loewner, data-based model reduction framework, as an extension of the Krylov methods for the case when the plant model is not available \cite{gosea2021datadriven,kergus-olivi-poussotVassal-demourant-LCSS2019}.

In this paper, we study the problem of designing a low order cloesed
loop system resulting in a low order controller achieving the closed-loop
stability and, guaranteeing asymptotic performance, assumed goals
given in a closed-loop desired model. In the proposed method,the control
design is based on . recasting the asymptotic performance problem
in a  closed-loop, time-domain moment matching based model reduction
framework. This problem has been discussed in terms of \emph{model}-matching
controllers for disturbance rejection \cite{medvedev-hillerstrom-CTAT1995}, where
the application of the Exact Model Matching concept, as in \cite{hautus-1980},
led to an algebraic solution to the Internal Model Principle, of which
more details can be found in, e.g., \cite{doyle-francis-tannenbaum-1992,francis-wonham-AUT1976,medvedev-hillerstrom-CTAT1995}
. Compared to \cite{medvedev-hillerstrom-CTAT1995}, we show that
the Internal Model Principle \textit{is equivalent to a time-domain
moment matching problem} and a stable parametrized reduced-order closed-loop
system, which achieves moment matching at the modes of the (persistent)
reference signal of the closed-loop, guarantees asymptotic tracking.After
designing a stabilizing reduced-order controller, based on moment
matching techniques, the asymptotic performance can be achieved through
a specific choice of the moments to match, i.e., the moments of the
closed-loop system match the moments of the signal generator. Stability
of the closed-loop system approximation is maintained through a proper
selection of the free parameters resulted from the application of
the moment matching technique. These results pave the way to obtain
low-order controllers based on the highly computationally efficient
Krylov projection-based moment matching techniques, as described in,
e.g., \cite{gugercin-antoulas-beattie-gildin-CDC2004}, to guarantee
closed-loop asymptotic performance.

The paper is organized as follows. In Section \ref{sec:Problem-formulation},
we present the problem formulation. In Section \ref{sect_prel}, we
give a brief overview of the recent results on time-domain moment
matching. In Section \ref{sec:Controller-reduction-track}, we present
our main results, i.e., asymptotic tracking is equivalent to achieving
moment matching at the interpolation points provided by the poles
of the signal generator. First we present the necessary and sufficient
condition for a linear system to asymptotically track a reference
signal and describe the family of low order models that achieve this
performance. Then we apply the results to the closed-loop  problem.
. In Section \ref{sect_example}, we illustrate the theory with an
example, followed by Discussion in Section \ref{sec:Discussion} and
Conclusions in Section \ref{sec:Conclusions}.

\paragraph*{Notation.}

$\mathbb{R}$ is the set of real numbers and $\mathbb{C}$ is the
set of complex numbers. $\mathbb{C}^{0}$ is the set of complex numbers
with zero real part and $\mathbb{C}^{-}$ denotes the set of complex
numbers with negative real part. $A^{*}$ denotes the conjugate transpose
of the complex matrix $A$. The transpose of $A$ is denoted by $A^{T}$.
$\sigma(A)$ denotes the set of eigenvalues of the matrix $A$ and
$\emptyset$ denotes the empty set. $$P=\left[\begin{array}{c|c}
A & B\\
\hline C & D
\end{array}\right]$$ denotes the transfer function $P(s)=C(sI-A)^{-1}B+D$. $\dim P$
is the dimension of a state-space realization of the transfer function
$P(s)$, i.e., the number of states.

\begin{figure}[H]
\label{fig_classic_CL}\centering
\includegraphics[width=.48\textwidth]{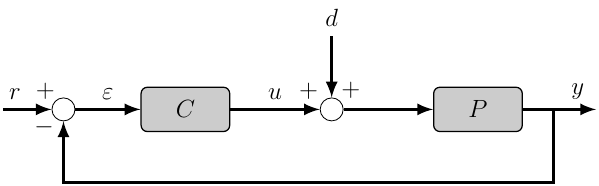}
\protect\caption{Classic feedback control system}
\end{figure}

\section{Problem formulation\label{sec:Problem-formulation}}

Consider the classical feedback control system shown in Figure \ref{fig_classic_CL},
where, for instance, a stabilizing controller $K$ has been designed
for a given plant $P$, i.e., we have a given closed-loop model. Typically,
the controller $K$ has an order, say $n_{C}$, that is much higher
than the plant order $n$, making them difficult to compute and implement.
 Let $\epsilon(t)=r(t)-y(t)$ be the error signal between the reference
and the output and $d(t)$ be a signal from a class of known disturbances.
The closed-loop transfer function $P_{CL}$ from the reference $r$
to the output $y$ can be described as 
\begin{equation}\label{eq:Pcl}
P_{CL}(s)=\frac{P(s)K(s)}{1+P(s)K(s)}.
\end{equation}
The problem to be studied in this paper is formulated below. 

\begin{prob} \label{prob_Kc_red} For the stable closed-loop
system shown in Figure \ref{fig_classic_CL}, find a family of $\nu<n_{c}$-th
order controllers $\widehat{K}$ and implicitly a family of low-order
closed-loop systems 
\begin{equation}\label{eq:PCL_hat}
\widehat{P}_{CL}(s)=\frac{P(s)\widehat{K}(s)}{1+P(s)\widehat{K}(s)}
\end{equation}
such that the closed-loop system $\widehat{P}_{CL}$ is asymptotically
stable, and achieves asymptotic reference tracking and disturbance
rejection, i.e., $\lim_{t\to\infty}\epsilon(t)=0$ for the given reference
$r(t)$ and the known type of disturbances $d(t)$. 

\end{prob} 

In this way, it is possible to decouple the stabilization problem
from achieving asymptotic performance, that is to design a low-order
stabilizing controller using moment matching, so that the closed-loop
system is (asymptotically) stable and satisfies other design features,
e.g., an optimal index, and then find a low-order approximation of
the designed controller such that the closed-loop system is (asymptotically)
stable and achieves desired asymptotic performance. 

The results are developed in the state-space framework, irrespective
of the computation of transfer functions.

\section{Preliminaries\label{sect_prel}}

In this section, the idea of moment matching for linear, single-input,
single-output (SISO) systems, from a time domain point of view as
presented in \cite{astolfi-CDC2010}, is outlined. Note that these
results are directly applicable to the multiple-input, multiple-output
case, see, e.g., \cite{i-astolfi-colaneri-SCL2014}.

\subsection{The notion of moment}

Consider the minimal linear, time-invariant SISO system 
\begin{equation}
P(s)=\left[\begin{array}{c|c}
A & B\\
\hline C & 0
\end{array}\right],\label{realization}
\end{equation}
with $\dim P=n$. 

\begin{deff}\n\label{def_moment}\cite{astolfi-TAC2010,antoulas-2005}
The 0-moment of system (\ref{realization}) at $s_{1}\in\mathbb{C}$
is the complex number $\eta_{0}(s_{1})=P(s_{1})$. The $k-$moment
of system (\ref{realization}) at $s_{1}$ is the complex number $${\displaystyle \eta_{k}(s_{1})=\left.\frac{(-1)^{k}}{k!}\cdot\frac{d^{k}P(s)}{ds^{k}}\right|_{s=s_{1}},}$$
with the integer $k\geq1$. \fin\end{deff}

Assume matrices $S\in\mathbb{R}^{\nu\times\nu}$, $L\in\mathbb{R}^{1\times\nu}$,
with $\nu\leq n$ and $\sigma(A)\cap\sigma(S)=\emptyset$, form an
observable pair $(L,\,S)$. Then, according to \cite{astolfi-TAC2010,astolfi-CDC2010,gallivan-vandendorpe-vandooren-SIAM2004,gallivan-vandendorpe-vandooren-JCAM2004},
the moments of system (\ref{realization}) at $\{s_{1},s_{2},...,s_{\nu}\}=\sigma(S)$
are given by $C\Pi$, where $\Pi\in\mathbb{C}^{n\times\nu}$ with
${\rm rank}\ \Pi=\nu$ is the unique solution of the Sylvester equation
\begin{equation}
A\Pi+BL=\Pi S,\label{eq_Sylvester}
\end{equation}
since the system \eqref{realization} is assumed minimal and the pair
$(L,S)$ is observable, as in, e.g., \cite{desouza-bhattacharyya-LAA1981}. 

\begin{obs}\n For the sake of clear exposition we have chosen the
matrix $S$ to have a spectrum of distinct values. However, this need
not be the case. The spectrum of $S$ may contain eigenvalues of higher
order multiplicities resultingin the presence of Jordan blocks in
the matrix $S$, see, e.g., \cite{astolfi-TAC2010,i-astolfi-colaneri-SCL2014}
for further details.\fin

\end{obs}

\subsection{Moment matching based model reduction\label{sub:Moment-matching}}

Let
\begin{equation}
\Sigma_{G}=\left[\begin{array}{c|c}
F & G\\
\hline H & 0
\end{array}\right]\label{sys_FGH}
\end{equation}
be a linear system with $\dim\Sigma_{G}=\nu\leq n$. Then its moments
match those of (\ref{realization}) at $\sigma(S)$ if 
\begin{equation}
HP=C\Pi,\label{MM_condition}
\end{equation}
where $P\in\mathbb{C}^{\nu\times\nu}$ is the unique solution of the
Sylvester equation 
\[
FP+GL=PS.
\]
Since the pair $(L,\,S)$ is assumed observable, a family of $\nu$th-order
models (\ref{sys_FGH}) that approximate (\ref{realization}) can
be described by 
\begin{equation}
\Sigma_{G}=\left[\begin{array}{c|c}
S-GL & G\\
\hline C\Pi & 0
\end{array}\right],\label{red_mod_CPi}
\end{equation}
parametrized in $G$, such that $\sigma(S-GL)\cap\sigma(S)=\emptyset$.

\subsection{Relation between moments  and the (well-defined) steady-state response}

The following result establishes the link between the moments of a
given system and its (well-defined) steady-state response when excited
by a given input signal. 

\begin{thm}\label{thm_MM_steady}\cite{astolfi-TAC2010} Consider
system (\ref{realization}) with $\sigma(A)\subset\mathbb{C}^{-}$.
Assume its input $u(t)$ is fed by the output $\theta(t)$ of the
signal generator 
\begin{align}\label{signal_gen}
\dot{\omega}&=S\omega,\\
\theta&=L\omega,\,\omega(0)\ne0, \nonumber
\end{align}
with $\omega(t)\in\mathbb{R}^{\nu}$ and $\theta\in\mathbb{R}$, where
the triplet $(L,S,\omega(0))$ is minimal and $\sigma(S)\subset\mathbb{C}^{0}$.
Then the moments $C\Pi$ of system (\ref{realization}) at $\sigma(S)$
are in a one-to-one relationship with the (well-defined) steady-state
response\footnote{By one-to-one relationship between a set of $\nu$ moments $\eta(s_{i})$,
$i=1,...,\nu$, and the well-defined steady state response of the
signal $y(t)$, we mean that the moments $\eta$ uniquely determine
the steady-state response of $y(t)$ , $y_{p}(t)=Cx_{p}(t)$, where
${\displaystyle x_{p}(t)=\lim_{t_{0}\to\infty}\int_{t_{0}}^{t}e^{A(t-\tau)}Bu(\tau)d\tau}$
is the permanent component of the state $x$.} of $y(t)$.\fin\end{thm} The value of $\nu$ is determined by the
size of the prescribed set of interpolation points. The notation and
the terminology are consistent with the equations used for exogenous
signal generators from \cite{francis-wonham-AUT1976}.

\section{Moment matching based reduction to achieve asymptotic performance\label{sec:Controller-reduction-track}}

\subsection{Open-loop reduction to achieve asymptotic tracking}

In this subsection, we consider the system (\ref{realization}) fed
by the signal generator (\ref{signal_gen}) and propose the condition
for the linear system (\ref{realization}) to asymptotically track
its input signal.

\begin{thm}\label{thm_CpiL} Consider the linear system (\ref{realization}),
with $\sigma(A)\subset\mathbb{C}^{-}$, fed by the signal generator
(\ref{signal_gen}) with the realization $(L,\,S,\,\omega(0))$ minimal,
i.e., $\theta(t)=u(t)$. Assume that the output $\theta(t)$ is a
persistent signal and $\sigma(A)\cap\sigma(S)=\emptyset$. Let $\epsilon(t)=y(t)-\theta(t)$.
Then $\lim_{t\to\infty}\epsilon(t)=0$ if and only if 
\[
L=C\Pi,
\]
where $\Pi$ is the unique solution of the Sylvester equation (\ref{eq_Sylvester}). \fin

\end{thm}

In other words, the steady-state response $y(t)$ of the system (\ref{realization})
equals to the steady-state response $\theta(t)$ of the signal generator
(\ref{signal_gen}) if and only if the moments of system (\ref{realization})
at $\sigma(S)$ are equal to the output map $L$ of the signal generator
(\ref{signal_gen}).

\begin{pf} The arguments from the proof of Theorem \ref{thm_MM_steady}
can be followed. The cascaded connection of the signal generator (\ref{signal_gen})
and the linear system (\ref{realization}), defined by $\theta(t)=u(t)$
, with the output $\epsilon(t)$ is described by
\begin{equation}
\begin{cases}
\dot{\omega} & \!\!\!\!=S\omega,\\
\dot{x} & \!\!\!\!=Ax+BL\omega,\\
\epsilon & \!\!\!\!=Cx-L\omega.
\end{cases}\label{eq:series_intercon}
\end{equation}
Since, by the assumption that $\sigma(A)\cap\sigma(S)=\emptyset$,
$\Pi$ is the unique solution of (\ref{eq_Sylvester}), then $x=\Pi\omega$
is a \foreignlanguage{british}{centre} manifold. Note that $\dot{x}-\Pi\dot{\omega}=A(x-\Pi\omega)$
is attractive since $\sigma(A)\subset\mathbb{C}^{-}$. Hence $\epsilon(t)=(C\Pi-L)\omega(t)+Ce^{At}(x(0)-\Pi\omega(0))$
and the claim follows. \fin

\end{pf}

The following result gives some significant tracking properties for
any system observing the conditions in Theorem \ref{thm_CpiL}. 

\begin{cor}\label{cor_track_multiple_sig}Let 
\begin{equation}
S=\begin{bmatrix}S_{1} & S_{3}\\
0 & S_{2}
\end{bmatrix},\,L=\begin{bmatrix}L_{1} & L_{2}\end{bmatrix},\label{eq:SL_block}
\end{equation}
with $S_{1}\in\mathbb{C}^{\nu_{1}\times\nu_{1}},\,S_{2}\in\mathbb{C}^{\nu_{2}\times\nu_{2}},\,L_{1}\in\mathbb{C}^{1\times\nu_{1}}$
and $L_{2}\in\mathbb{C}^{1\times\nu_{2}}$. Assume the pair $(L,\,S)$
is observable, with $\dim S=\nu=\nu_{1}+\nu_{2}$. Consider the cascaded
connection of the signal generator (\ref{signal_gen}) and the system
(\ref{realization}), with $\sigma(A)\subset\mathbb{C}^{-}$, defined
by $\theta(t)=u(t)$. If the pairs $(L_{1},S_{1})$ and $(L_{2},S_{2})$
are observable, respectively, then the following statements are equivalent:
\begin{enumerate}
\item $L=C\Pi$, where $\Pi=[\Pi_{1}\,\,\Pi_{2}]$ is the unique solution
of the Sylvester equation (\ref{eq_Sylvester}).
\item $L_{1}=C\Pi_{1}$ and $L_{2}=C\Pi_{2}$.
\item System (\ref{realization}) asymptotically tracks the output signal
$\theta(t)$ of the signal generator (\ref{signal_gen}).
\item System (\ref{realization}) asymptotically tracks each of the signals
of the signal generators 
\begin{align*}
\dot{\omega}_{1} & =S_{1}\omega_{1},\ \theta_{1}=L_{1}\omega_{1},\\
\dot{\omega}_{2} & =S_{2}\omega_{2},\ \theta_{2}=L_{2}\omega_{2},
\end{align*}
respectively.\fin
\end{enumerate}
\end{cor}

\begin{pf} The proof is a direct consequence of Theorem \ref{thm_CpiL},
using the block-version solution of the Sylvester equation (\ref{eq_Sylvester}).
\fin

\end{pf}

In particular, as a direct consequence of Corollary \ref{cor_track_multiple_sig},
the following result holds for polynomial signals $t^{k}$, $k\in\mathbb{N}$.

\begin{cor}\label{cor_tn} If the system (\ref{realization}), with
$\sigma(A)\subset\mathbb{C}^{-}$, asymptotically tracks the signal
$t^{\nu},\,t\geq0$, then the system (\ref{realization}) asymptotically
tracks each of the signals $t^{k},\,t\geq0$, for all $k=0,\dots,\nu$.
\fin

\end{cor}

Finally, consider the subfamily of asymptotically stable models $\Sigma_{G}$,
\[
\bar{\Sigma}_{G}=\left\{ \Sigma_{G}\,{\rm as\,in\,(\ref{red_mod_CPi})\,}|\,G\,{\rm is\,such\,that\,}\sigma(S-GL)\subset\mathbb{C}^{-}\right\} .
\]
The following result yields a condition such that a system $\bar{\Sigma}_{G}$
achieves asymptotic tracking. 

\begin{prop}\label{prop_GL} For $H=L$, where $H$ is as in (\ref{sys_FGH}),$\bar{\Sigma}_{G}$
defines a family of systems of order $\nu$ that asymptotically track
the output signal $\theta(t)$ of the signal generator (\ref{signal_gen}).
Furthermore, the models $\bar{\Sigma}_{G}$ match the moments of (\ref{realization})
at $\sigma(S)$ if there exists an invertible matrix $\widetilde{P}$
such that $C\Pi=L\widetilde{P}$.\fin

\end{prop} 

\begin{pf}The proof follows the arguments similar to the proof of
Theorem \ref{thm_CpiL}. Consider a signal generator described by
the pair $(C\Pi,\,S)$. The final statement follows from the moment
matching condition, i.e., (\ref{realization}) matches the moments
of $\bar{\Sigma}_{G}$ at $\sigma(S)$ if there exists an invertible
matrix $\widetilde{P}$ such that $L\widetilde{P}=C\Pi.$ Note that
since $(L,\,S)$ is observable, the matrix $\widetilde{P}$ exists
and the claim follows.\fin

\end{pf}

\subsection{Closed-loop reduction to achieve stabilization and asymptotic performance\label{sub:Closed-loop-reduction}}

In this subsection, we consider the classical feedback control system
shown in Figure \ref{fig_classic_CL}, where the reference signal
$r(t)$ is yielded by the signal generator (\ref{signal_gen}). Furthermore,
the error signal is defined as $\epsilon(t)=r(t)-y(t)$. The closed-loop
transfer function from $d(t)$ to $y(t)$ is
\begin{equation}
T(s)=\frac{P(s)}{1+P(s)K(s)},\label{eq:Ts}
\end{equation}
for the control scheme in Figure \ref{fig_classic_CL}, with the controller
$K$ and 
\begin{equation}
\widehat{T}(s)=\frac{P(s)}{1+P(s)\widehat{K}(s)},\label{eq:Ts_hat}
\end{equation}
for the control scheme in Figure \ref{fig_classic_CL}, with the controller
$\hat{K}$. The transfer function from the reference $r(t)$ to the
error $\epsilon(t)$ is
\begin{equation}
E(s)=\frac{1}{1+P(s)K(s)}.\label{eq:Es}
\end{equation}
for the control scheme in Figure \ref{fig_classic_CL}, with the controller
$K$ and 
\begin{equation}
\widehat{E}(s)=\frac{1}{1+P(s)\widehat{K}(s)},\label{eq:Es_hat}
\end{equation}
for the control scheme in Figure \ref{fig_classic_CL}, with the controller
$\hat{K}$.  A few properties are in order at this point. 

\begin{lemma}\label{lema_interp} Let $s_{0}\in\mathbb{C}$. If $K(s_{0})=\widehat{K}(s_{0})$
then $P_{CL}(s_{0})=\widehat{P}_{CL}(s_{0})$, $T(s_{0})=\widehat{T}(s_{0})$
and $E(s_{0})=\widehat{E}(s_{0})$, where $\widehat{P},\,\widehat{T}$
and $\widehat{E}$ are as in (\ref{eq:PCL_hat}), (\ref{eq:Ts_hat})
and (\ref{eq:Es_hat}), respectively. \fin

\end{lemma}

\begin{pf}  The results is rather straightforward, hence the proof
is omitted.\fin

\end{pf} If $s_{0}$ is not a zero of the plant, i.e., $P(s_{0})\ne0$,
then the relation $K(s_{0})=\widehat{K}(s_{0})$ in Lemma \ref{lema_interp}
becomes a necessary and sufficient condition.

The result of Theorem \ref{thm_CpiL} is applicable to the closed-loop
case, as follows.

\begin{thm}\label{thm_KclT} Let $\widehat{P}_{CL}$ as in (\ref{eq:Pcl})
and $\widehat{E}(s)$ as in (\ref{eq:Es}) be the closed-loop transfer
function and the sensitivity transfer function of the negative feedback
control loop shown in Figure \ref{fig_classic_CL}, where $K$ is
replaced by $\widehat{K}$. Assume that $\widehat{P}_{CL}$ and implicitly
$\widehat{E}$ are asymptotically stable. Let the reference $r(t$)
be the output of the signal generator (\ref{signal_gen}) with the
pair $(L,\,S)$ observable. If $\dim\widehat{P}_{CL}=\dim\widehat{E}=\dim S$,
then the following statements are equivalent:
\begin{enumerate}
\item $\widehat{P}_{CL}$ asymptotically tracks the reference $r(t)$.
\item The moments of $\widehat{P}_{CL}$ at $\sigma(S)$ are equal to the
elements of $L$.
\item The moments of $\widehat{E}(s)$ at $\sigma(S)$ are equal to zero.
\fin 
\end{enumerate}
\end{thm}

\begin{pf} $(1)\Leftrightarrow(2)$ follows immediately from Lemma
\ref{lema_interp} and Theorem \ref{thm_CpiL}. $(2)\Leftrightarrow(3)$
follows from a Final Value Theorem argument and using the fact that
$E(s)=1-P_{CL}(s)$ and the claim follows.\fin

\end{pf}

Similar arguments hold for $T(s)$ as in (\ref{eq:Ts}), proving identical
results for asymptotic disturbance rejection. Note that Theorem \ref{thm_KclT}
is equivalent to the usual Internal Model Principle, as in e.g., \cite{doyle-francis-tannenbaum-1992}, since,
in order to achieve asymptotic tracking of the reference signal, the
model of the reference ought to be included in the closed-loop system.
In our case, we actually start the reduced order design from achieving
moment matching at the interpolation points provided by the reference
signal generator. These results pave the way for the design of low-order
controllers, using efficient Krylov projection-based algorithms, as
in, e.g., \cite{gugercin-antoulas-beattie-gildin-CDC2004}, with a
specific choice of interpolation points given by the poles of the
reference signal generator.

A closed-loop reduction algorithm may be described as follows. Consider
the plant $P$, as in (\ref{realization}) and let $(L,S)$ be as
in (\ref{eq:SL_block}). The dimension of the reduced closed-loop
system is determined by the selection of the dimension of the signal
generator used in the moment matching procedure, $\dim S=\nu=n+\nu_{C}$,
where $\nu_{C}$ is the dimension of the reduced order controller,
i.e.,$\dim S_{2}=\nu_{C}$. 

For a controller realized as
\[
K=\left[\begin{array}{c|c}
A_{K} & B_{K}\\
\hline C_{K} & D_{K}
\end{array}\right],
\]
 the transfer function of the closed-loop system is 
\[
P_{CL}=\left[\begin{array}{c|c}
A_{CL} & B_{CL}\\
\hline C_{CL} & D_{CL}
\end{array}\right]=\left[\begin{array}{cc|c}
A-BD_{K}C & BC_{K} & BD_{K}\\
-B_{K}C & A_{K} & B_{K}\\
\hline C & 0 & D
\end{array}\right],
\]
with $\dim P_{CL}=n+n_{c}$. Applying the moment matching ideas described
in Section \ref{sub:Moment-matching} yields an approximation $\hat{P}_{CL}$,
with $\dim\widehat{P}_{CL}=n+\nu_{C}$, described by
\begin{equation}
\widehat{P}_{CL}\!\!=\!\!\left[\begin{array}{c|c}
\widehat{A}_{CL} & \widehat{B}_{CL}\\
\hline \widehat{C}_{CL} & \widehat{D}_{CL}
\end{array}\right]\!\!=\!\!\left[\begin{array}{cc|c}
S_{1}-G_{1}L_{1} & G_{1}C_{K}\Pi_{2} & G_{1}\\
S_{3}-G_{2}C\Pi_{1} & S_{2}-G_{2}L_{2} & G_{2}\\
\hline C\Pi_{1} & 0 & D
\end{array}\right],\label{eq:Pcl_red}
\end{equation}
where $G=[G_{1}^{*}\,\,G_{2}^{*}]^{*}\in\mathbb{C}^{n+\nu}$ are free
parameters and $\Pi_{1}$ is such that $(A-BD_{K}C)\Pi_{1}=\Pi_{1}S_{1}$
and $\Pi_{2}$ is the unique solution of the Sylvester equation (\ref{eq_Sylvester})
$A_{K}\Pi_{2}+B_{K}L_{2}=\Pi_{2}S_{2}$. Letting $G_{1}=B$, $\widehat{P}_{CL}$
is the closed-loop system formed of the given plant with the reduced
order controller. Furthermore, setting $C\Pi_{1}=L_{1}$, by Theorem
\ref{thm_CpiL}, if $\bar{P}_{CL}$ matches the moments at the reference
signal generator, characterized by $(C\Pi_{1},S_{1})$, then the reference
is asymptotically tracked.

\subsection{A closer look at the condition $C\Pi=L$.}

In this subsection, if the condition $C\Pi=L$ causes any loss of
generality is investigated. If $L=C\Pi$, the Sylvester equation (\ref{eq_Sylvester})
becomes 
\begin{equation}
(A+BC)\Pi=\Pi S.\label{eq:Sylv_CPi}
\end{equation}
This equation has non-zero solutions if 
\begin{equation}
\sigma(A+BC)\cap\sigma(S)\ne\emptyset.\label{eq:ker}
\end{equation}
Note that this condition does not restrict the class of linear systems
that satisfy Theorem \ref{thm_CpiL}, since the restriction imposed
by (\ref{eq:ker}) is on the choice of $S$ ,i.e., the choice of the
signal generator is slightly restricted. Furthermore, equation (\ref{eq:Sylv_CPi})
can be recast using the Kronecker product as (see, e.g., \cite{golub-nash-vanloan-TAC1979,hu-LAA1992})
\[
(I_{n}\otimes(A+BC)-S\otimes I_{\nu}){\rm vec}(\Pi)=0,
\]
where ${\rm vec}(\Pi)=[\Pi_{1}^{T}\,\Pi_{2}^{T}\,\dots\,\Pi_{\nu}^{T}]^{T}\in\mathbb{C}^{\nu n}$,
with $\Pi_{k}\in\mathbb{C}^{n},\,k=1,\dots\nu$ is a column of matrix
$\Pi$. Hence, any solution $\Pi$ of (\ref{eq:Sylv_CPi}) satisfies
the condition 
\[
{\rm vec}(\Pi)\in\ker(I_{n}\otimes(A+BC)-S\otimes I_{\nu}).
\]
Note that, indeed a zero unique solution exists if and only if the
matrix $(I_{n}\otimes(A+BC)-S\otimes I_{\nu})$ is invertible, i.e.,
$\sigma(A+BC)\cap\sigma(-S)=\emptyset.$ 

In the case of system (\ref{sys_FGH}), condition (\ref{eq:Sylv_CPi})
becomes $(F+GH)P=PS$. For $F=S-GL$, condition (\ref{eq:ker}) becomes
$\sigma(F+GH)=\sigma(S)$. Hence, choosing $H=L$, generally defines
a family of models $\Sigma_{G}$ that asymptotically tracks signal
$\theta(t)$ yielded by the signal generator (\ref{signal_gen}).

\section{Illustrative example}

\label{sect_example}

\begin{figure}
\includegraphics[clip]{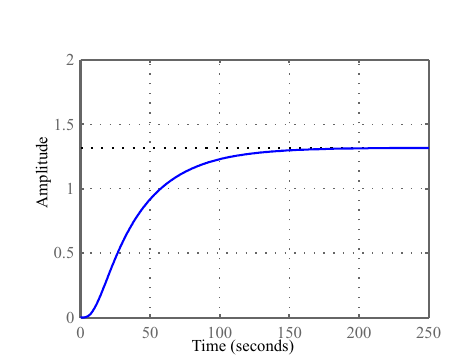} \protect\caption{Step response of an 8th-order four disk drive system with a Kalman
compensator. Asymptotic tracking is not achieved.}

\label{CP_t2-1}
\end{figure}

The four-disk drive system studied in \cite{anderson-liu-ACC1987,liu-anderson-IJC1987}
is taken as an example. The plant is described by an unstable LTI
SISO system as in (\ref{realization}) of eighth order with{\tiny
\begin{multline*}
A=\begin{bmatrix}-0.161 & -6.004 & -0.58215 & -9.9835 & -0.40727 & -3.982 & 0 & 0\\
1 & 0 & 0 & 0 & 0 & 0 & 0 & 0\\
0 & 1 & 0 & 0 & 0 & 0 & 0 & 0\\
0 & 0 & 1 & 0 & 0 & 0 & 0 & 0\\
0 & 0 & 0 & 1 & 0 & 0 & 0 & 0\\
0 & 0 & 0 & 0 & 1 & 0 & 0 & 0\\
0 & 0 & 0 & 0 & 0 & 1 & 0 & 0\\
0 & 0 & 0 & 0 & 0 & 0 & 1 & 0
\end{bmatrix},\\
B^{T}=[1\,\,0\,\,0\,\,0\,\,0\,\,0\,\,0\,\,0],\\
C=[0\,0\,6.4432\cdot10^{-3}\,2.1936\cdot10^{-3}\,7.1252\cdot10^{-2}\,1.0002\,0.10455\,0.99551].
\end{multline*}
}In this example, we compare the propsed moment matching-based closed-loop
design with a classic multivariable compensator, e.g. given in \cite[Chapter 7, Section 7.5]{astrom-murray-2012},
which yields a closed-loop system of order 16. The eighth order Kalman
compensator stabilizes the closed-loop system with the eigenvalues
-1.5000, -1.1000, -1.0000, -0.5000, -0.3333, -0.2500, -0.2000, -0.1000,
-2.1000, -1.3000, -1.0000, -0.2000, -0.1667,-0.1429, -0.0300, -0.0100.
With this controller, the closed-loop system does not achieve asymptotic
performance, as shown in Figure \ref{CP_t2-1}. 

\begin{figure}
\includegraphics{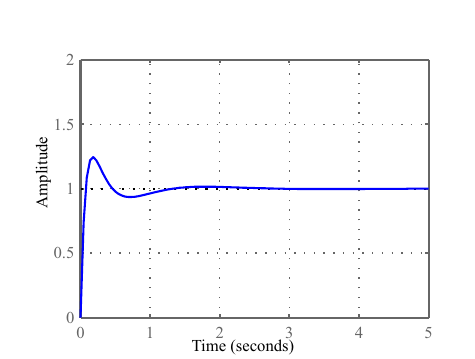}\protect\caption{Step response of an 8th-order four disk drive system with a fourth
order controller. Asymptotic tracking of a step reference signal is
achieved.}

\label{CP_step-1}
\end{figure}

We choose $\nu_{C}=4$ based on empirical observation. Hence $\nu=n+\nu_{C}=12$,
i.e., the first twelve moments at 0 are chosen. Let $S$ be as in
\ref{eq:SL_block}, where $S_{1}$ is a Jordan block of order 8 with
0 an eigenvalue of multiplicity 8 and $S_{2}$ is a Jordan block of
order 4 with 0 an eigenvalue of multiplicity 4. By the algorithm given
in Section \ref{sub:Closed-loop-reduction}, a fourth-order controller
is obtained. The family of twelve-order closed-loop reduced order
systems $\widehat{P}_{CL}$ are described by (\ref{eq:Pcl_red}),
with a fourth order controller. We have computed $G$ and found a
stable $\bar{P}_{CL}$, with the closed-loop eigenvalues $-8.6072$,
$-0.0370+1.8496j$, $-0.0370-1.8496j$, $-0.0282+1.4097j$, $-0.0282-1.4097j$,
$-1.3845$, $-0.0153+0.7648j$, $-0.0153-0.7648j$, $-0.02$, $-0.001$,
$-0.0042+0.0082j$, $-0.0042-0.0082j$. Based on Theorem \ref{thm_KclT},
the closed-loop system $\bar{P}_{CL}$, with the fourth order reduced
controller, asymptotically tracks a step reference signal. This is
demonstrated by the step response shown in Figure \ref{CP_step-1}.

\section{Discussion\label{sec:Discussion}}

A few remarks are in order.

We apply recent result on time-domain moment matching to controllers
in closed-loop systems as it has not been applied before.

We connect moment-matching based model/controller reduction to the
Internal Model Principle, loweing the order of the controller but
maintaining the stability of the closed-loop system and, the asymptotic
performance.

The Sylvester equation is only used for the conceptual development
of the proposed method. The solution need not be computed at all,
as it can be obtained using Krylov projection methods, see, e.g.,
\cite{i-astolfi-colaneri-SCL2014,astolfi-CDC2010}, for more details.
Hence the results may be applied to larger dimensional systems.

The moment matching-based approach exhibits additional degrees of
freedom proven to be useful for yielding suitable reduced-order systems.
We do not provide any systematic way of selecting the degrees of freedom
in all cases, nor do we provide a theoretical approach on how to select
an appropriate approximation order. There is no precedent on these
issues in the LTI case whatsoever. We limit ourselves to selecting
the free parameters on a trial and error and experience basis. 

We do not make any theoretical analysis on the accuracy of the approximation,
which is intrinsic to moment matching-based approximations. However,
the degrees of freedom allow for choices of interpolation points and
other free parameters such that lower errors are attained. Again,
this is done on simulation basis. Similar arguments apply to the choice
of interpolation points used in the selection of the signal generator
pair $(L,S)$. Furthermore, the dimension of the reduced closed-loop
system and of the lower order controller is dictated by the number
of interpolation points selected. To our knowledge, currently, there
exists no systematic way of determining the dimension of the low order
approximant. It is usually done based on experience. However, recently,
in \cite{ionescu-scherpen-iftime-astolfi-MTNS2012}, preliminary results
establish a relation of the moments to specific Hankel singular values,
in view of solving the problem of appropriate choices of reduced order
dimensions.

\section{\label{sec:Conclusions}Conclusions}

In this paper we have studied the problem of model reduction of closed-loop
systems with the goal of achieving desired asymptotic performance.
We have established conditions such that moment matching of the output
map of the reference signal generator is identical to finding a low
order closed-loop system that achieves asymptotic tracking of the
reference. These results pave the way for design of low order controllers,
using efficient Krylov-projection based algorithms with a specific
choice of interpolation points given by the poles of the reference
signal generator.

\bibliographystyle{plain}
\bibliography{BIB/TCI_articles,BIB/TCI_books,BIB/TCI_phdthesis}

\end{document}